\documentclass[12pt,a4paper]{article}
\usepackage{theorem}
\newtheorem{lemma}{Lemma}
\newtheorem{proposition}[lemma]{Proposition}
\newtheorem{theorem}[lemma]{Theorem}

{\theorembodyfont{\upshape}}
{\theorembodyfont{\upshape}}
{\theorembodyfont{\upshape}\newtheorem{example}[lemma]{Example}}
{\theorembodyfont{\upshape}}

\newcommand{\R}{{\bf R}}
\newcommand{\C}{{\bf C}}
\newcommand{\rme}{{\rm e}}
\newcommand{\rmd}{{\rm d}}

\newcommand{\cC}{{\cal C}}

\newcommand{\cH}{{\cal H}}

\newcommand{\sig}{\sigma}
\newcommand{\alp}{\alpha}
\newcommand{\bet}{\beta}
\newcommand{\gam}{\gamma}
\newcommand{\lam}{\lambda}
\newcommand{\del}{\delta}

\newcommand{\lap}{{\Delta}}
\newcommand{\gra}{\nabla}

\newcommand{\Dom}{{\rm Dom}}

\newcommand{\Spec}{{\rm Spec}}

\newcommand{\Ran}{{\rm Ran}}

\newcommand{\norm}{\Vert}
\renewcommand{\Re}{{\rm Re}\;}

\newcommand{\supp}{{\rm supp}}
\newcommand{\Proof}{\underbar{Proof}{\hskip 0.1in}}

\newcommand{\Num}{{\rm Num}}
\newcommand{\conv}{{\rm conv}}
\newcommand{\dist}{{\rm dist}}
\newcommand{\Schrodinger}{Schr\"odinger }

\newcommand{\abs}[1]{\left\vert#1\right\vert}

\newcommand{\la}{{\langle}}
\newcommand{\ra}{{\rangle}}

\title{THE PSEUDOSPECTRAL PROPERTIES OF NON-SELF-ADJOINT SCHR\"ODINGER
OPERATORS IN THE SEMI-CLASSICAL LIMIT}
\author{Paul Redparth}
\date{July 2000}
\begin{document}
\maketitle
\begin{abstract}
We describe the general qualitative behaviour of the resolvent
norm for a very wide class of non-self-adjoint Schr\"odinger
operators in the semi-classical regime, as the spectral parameter
$\lam$ varies over the complex plane.

\vskip 0.1in AMS subject classifications: 47A10, 81Q20.
\par
keywords: \Schrodinger operator, non-self-adjoint operator,
pseudospectrum, WKB approximation, semi-classical analysis.
\end{abstract}

\section{Introduction}\label{s1}

Several authors \cite{Dav97a,Reddy93,Tre97} have demonstrated that
the spectrum of certain Schr\"odinger operators with complex
potentials is unstable; this is shown by observing that the
resolvent norm $\norm(H-\lam)^{-1}\norm$ becomes unbounded as some
parameter associated with the operator $H$ varies, even though
$\lam$ may be far from the spectrum of $H$. This phenomena is
commonly quantified using the concept of the pseudospectral sets
\[
\Spec_{\epsilon}(H):=\Spec(H)\cup\{\lam\in\C:\norm(H-\lam)^{-1}\norm
\geq\epsilon^{-1}\},
\]
where we adopt the convention that if $\lam\in\Spec(H)$ then $
\norm(H-\lam)^{-1}\norm:=\infty.$ For any closed operator $T$
acting on a Hilbert space $\cH$ the $numerical$ $range$, defined
by
\[
\Num(T):=\{\la Tf,f\ra:f\in\Dom(T),\mbox{ and } \norm f\norm=1\}
\]
is a convex subset of $\C$ \cite[Theorem 6.1]{Dav80}, containing
$\Spec(T)$ if $(T+s)$ is $maximal$ $quasi$-$accretive$
\cite[p.279]{Kato66}; equivalently, if
\[
\Re(\Num(T+s))\geq 0\qquad\mbox{ and }\qquad\Ran(T+t)=\cH
\]
for some (and hence all) $t>s$. It is then well-known that,
provided $\C\backslash\overline{\Num(T)}$ is a connected set
\begin{equation}\label{res1}
\norm(T-\lam)^{-1}\norm\leq\frac{1}{\dist(\lam,\overline{\Num(T)})}
\end{equation}
for all $\lam\notin\Num(T)$ \cite[p.268]{Kato66}. We will be
concerned with the Schr\"odinger operator
\[
H_h:=-h^2\lap+V
\]
with complex-valued $V$ and $h>0$, acting in $L^2(\Omega)$ where
$\Omega$ is some region in $\R^N$. We assume Dirichlet boundary
conditions throughout. The $h$-independent set
\begin{equation}\label{def1}
\Phi(V):=\{\overline{\Ran(V)}+[0,\infty)\}
\end{equation}
will be of fundamental importance. In Section \ref{s2} we show
that for any $\lam\in\Phi(V)$, the resolvent norm tends to
infinity in the semi-classical limit, $h\to 0$. The fact that $V$
may take complex values means that
$\conv(\Phi(V))\backslash\Phi(V)$ (conv denoting the convex hull)
is in general non-empty, and in Section \ref{s3} we give new
results (Theorems \ref{t6} and \ref{t7}) which show that in the
semi-classical limit a bound analogous to (\ref{res1}) holds for
$\lam$ in this set.

\section{Preliminaries}\label{s1a}

First we discuss the conditions which we shall impose on $V$ and
the question of the domain of the operator $H_h$. We assume that
there exists a closed set $M\subseteq\Omega$ with Lebesgue measure
zero, such that the restriction $V\Big|_{\Omega\backslash M}$ is
continuous and bounded. Since the operator $H_h$ will not be
changed on sets of measure zero, we allow $V$ to be undefined on
$M$. Defining $H_{h,0}$ to be the self-adjoint operator $-h^2\lap$
with domain $W_0^{1,2}(\Omega)$, it is trivial that $V$ has
relative bound zero with respect to $H_{h,0}$ since $V\in
L^{\infty}(\Omega)$. Therefore, $\Dom(H_h)=W_0^{1,2}(\Omega)$, and
the fact that
\[
\norm V(H_{h,0}+\lam)^{-1}\norm\leq\norm V\norm
_{\infty}\norm(H_{h,0}+\lam)^{-1}\norm =\lam^{-1}\norm V\norm
_{\infty}
\]
for all $\lam >0$, together with an argument similar to the proof
of \cite[Theorem 1.4.2]{Dav95} shows that $H_h$ is maximal on
$W_0^{1,2}(\Omega)$. Moreover, since we may integrate by parts
\begin{equation}\label{diffparts}
\la H_hf,f\ra=\int_{\Omega}V(x)\abs{f(x)}^2\rmd
x+h^2\int_{\Omega}\abs{\gra f}^2\rmd x
\end{equation}
for all $f\in W_0^{1,2}(\Omega)$. Since
\[
\{\Re(V(x)):x\in\Omega\backslash M\}\geq k
\]
for some $k\in\R$, it follows that $H_h$ is maximal
$quasi$-$accretive$ and so
$\Spec(H_h)\subseteq\overline{\Num(H_h)}$. It also follows from
(\ref{diffparts}) that
\begin{lemma}\label{t3}
With $H_h$ as defined above
\[
\overline{\Num(H_h)}\subseteq\conv(\Phi(V))
\]
for all $h>0$.
\end{lemma}

\begin{proposition}\label{p1}
For all $\lam\notin\conv(\Phi(V))$
\[
\norm(H_h-\lam)^{-1}\norm\leq\frac{1}{\dist(\lam,\conv(\Phi(V)))}.
\]
\end{proposition}
\Proof $H_h$ satisfies the conditions for (\ref{res1}) to hold,
and then one can apply the lemma.

\section{Spectral instability}\label{s2}

In a sense our result in this section generalises that of
\cite[Section 2]{Dav98b}; however, our conclusion is not as strong
since we do not show super-polynomial growth of the resolvent norm
in the semi-classical limit. Note that the specific form of the
phase function $x\mapsto\gam\cdot x$ below has been chosen to
ensure that the second derivative disappears, and to exploit the
fact that $\lam:=V(c)+\abs{\gam}^2$. We have not attempted to
optimise the choice of phase function, but instead have aimed for
a clear demonstration of the underlying processes.
\begin{theorem}\label{t2}
For $z\in\Phi(V)$, we have
\[
\norm(H_h-\lam)^{-1}\norm\to\infty\qquad\mbox{ as }h\to 0.
\]
\end{theorem}
\Proof Let $\lam\in\{\Ran(V)+[0,\infty)\}$ so that
$\lam:=V(c)+\abs{\gam}^2$, where $c\in\Omega\backslash M$, and
$\gam\in\R^N$. Let $\phi\in C^{\infty}_c(\R^N)$ be a function
whose support is contained within the open ball $B(0;1)$, and put
\[
\phi_{\lam,h}(x):=\phi\left(\frac{x-c}{h^{p}}\right)
\]
where $p>0$ is a constant to be determined, and
\[
f_{\lam,h}(x):=\rme^{ih^{-1}\gam\cdot x}\phi_{\lam,h}(x).
\]
Then $\supp(f_{\lam,h})\subseteq B(c;h^{p})$ and
\begin{eqnarray*}
&&H_hf_{\lam,h}\\ &&=-h^2\lap\rme^{ih^{-1}\gam\cdot
x}\phi_{\lam,h} +V\rme^{ih^{-1}\gam\cdot x}\phi_{\lam,h}\\
                  &&=-h^2\rme^{ih^{-1}\gam\cdot x}\lap\phi_{\lam,h}
-h^2\phi_{\lam,h}\lap\rme^{ih^{-1}\gam\cdot x}-2h^2(\gra\rme^
{ih^{-1}\gam\cdot
x}\cdot\gra\phi_{\lam,h})+V\rme^{ih^{-1}\gam\cdot x}
\phi_{\lam,h}\\
                  &&=-h^2\rme^{ih^{-1}\gam\cdot x}\lap\phi_{\lam,h}
+\abs{\gam}^2\rme^{ih^{-1}\gam\cdot x}\phi_{\lam,h}-2ih
\rme^{ih^{-1}\gam\cdot x}\gam\cdot\gra\phi_{\lam,h}
+V\rme^{ih^{-1}\gam\cdot x}\phi_{\lam,h}.\\
\end{eqnarray*}
Therefore, it follows that
\begin{eqnarray*}
\norm(H_h-\lam)f_{\lam,h}\norm_2 &\leq& h^2
\norm\lap\phi_{\lam,h}\norm_2+
2h\norm\gam\cdot\gra\phi_{\lam,h}\norm_2
+\norm(V(x)-V(c))\Big|_{B(c;h^{p})}\norm_{\infty}\norm
\phi_{\lam,h}\norm_2\\
\end{eqnarray*}
and so
\begin{eqnarray*}
\frac{\norm(H_h-\lam)f_{\lam,h}\norm_2}{\norm f_{\lam,h}\norm_2}
&\leq &
\frac{h^2\norm\lap\phi_{\lam,h}\norm_2}{\norm\phi_{\lam,h}\norm_2}+
\frac{2h\norm\gam\cdot\gra\phi_{\lam,h}\norm_2}{\norm\phi_{\lam,h}\norm_2}
+\norm(V(x)-V(c))\Big|_{B(c;h^{p})}\norm_{\infty}.\\
\end{eqnarray*}
Taking the limit,
\begin{eqnarray*}
\lim_{h\to 0}\frac{\norm(H_h-\lam)f_{\lam,h}\norm_2} {\norm
f_{\lam,h}\norm_2} &\leq &\lim_{h\to
0}\left(k_1h^{2-2p}+k_2h^{1-p}+
\norm(V(x)-V(c))\Big|_{B(c;h^{p})}\norm_{\infty}\right)\\
\end{eqnarray*} where $k_1,k_2$ are constants dependent upon
our choice of $\phi$ and $\gam$.\newline Thus, by taking $0<p<1$,
the continuity of $V$ on $\Omega\backslash M$ ensures that
\[
\lim_{h\to 0}\frac{\norm(H_h-\lam)f_{\lam,h}\norm_2} {\norm
f_{\lam,h}\norm_2}=0,
\]
or equivalently, that
\[
\lim_{h\to 0}\norm(H_h-\lam)^{-1}\norm =\infty.
\]
Now let $w\in\Phi(V)$ and, aiming for a contradiction, suppose
that
\[
\norm(H_h-w)^{-1}\norm\leq m\qquad\mbox{ as }h\to 0.
\]
For any $\del>0$ there exists $\lam\in\{\Ran(V)+[0,\infty)\}$ such
that $\abs{\lam-w}<\del$, by the definition of the closure. Hence,
using the resolvent identity,
\begin{eqnarray*}
\norm(H_h-\lam)^{-1}\norm
&=&\norm(H_h-w)^{-1}-(\lam-w)(H_h-\lam)^{-1}(H_h-w)^{-1}\norm\\
&\leq& m+\abs{\lam-w}m\norm(H_h-\lam)^{-1}\norm\\ &<& m+\del
m\norm(H_h-\lam)^{-1}\norm\qquad\mbox{ as }h\to 0.\\\
\end{eqnarray*}
Since $\del$ may be taken arbitrarily small, we have
\[
\norm(H_h-\lam)^{-1}\norm\leq m
\]
as $h\to 0$, contradicting the result just obtained, and
completing the proof.

\section{Bounded behaviour of the Resolvent norm}\label{s3}

We have seen that when $\lam\in\Phi(V)$ the resolvent norm becomes
infinitely large as $h\to 0$ (Theorem \ref{t2}). Conversely, when
$\lam\notin\conv(\Phi(V))$ the resolvent norm is uniformly bounded
in $h>0$ (Proposition \ref{p1}). Therefore, the natural question
arises: how does $\lim_{h\to 0}\norm(H_h-\lam)^{-1}\norm$ behave
for $\lam\in\conv(\Phi(V))\backslash\Phi(V)$? On the one hand,
Proposition \ref{p1} cannot in general be extended to
$\lam\in\conv(\Phi(V))\backslash \Phi(V)$, as the following
counter-example shows:
\begin{example}\label{t5}
Consider the operator
\[
H_{\del,h}f(x):=-h^2 f''(x)+V_{\del}(x)f(x)\qquad\mbox{ acting on
}L^2(-1,1)
\]
where
\[
V_{\del}(x):=\left\{\begin{array}{ll} i(x+\del)&\mbox{ for $x>0$
}\\ i(x-\del)&\mbox{ for $x<0.$ }\\ \end{array} \right.
\]
For any $\lam>0$, it follows that
$\lam\in\conv(\Phi(V_{\del}))\backslash\Phi(V_{\del}),$ where
\[
\Phi(V_{\del}):=\{\overline{\Ran(V_{\del})}+[0,\infty)\}.
\]
But, we have shown elsewhere \cite{Red1} that countably many of
the eigenvalues $\{\lam_{h,n}\}_{n=1}^{\infty}$ of $H_{\del,h}$
are positive real, for $every$ $h>0$, in which case
\[
\norm(H_{\del,h}-\lam_{h,n})^{-1}\norm :=\infty.
\]
\end{example}

On the other hand, for a wide class of potentials, we have the
positive result (initially suggested by numerical simulations of
the associated discrete problem using Matlab) that for any given
$\lam\in\conv(\Phi(V))\backslash\Phi(V)$ the resolvent norm
becomes bounded eventually, as $h\to 0$. In the one-dimensional
case we have the following result, which we believe to be new.

\begin{theorem}\label{t6}
Let $K_h$ be the non-self-adjoint operator
\begin{equation}\label{K}
K_hf(x):=-h^2\frac{\rmd^2f}{\rmd x^2}+V(x)f(x)
\end{equation}
acting in $L^2(a,b)$, $-\infty\leq a<b\leq +\infty$ and $h>0$.
When $a$ or $b$ are finite we impose Dirichlet boundary
conditions. We assume that there exists a partition
\[
a=x_0<x_1<\cdots<x_n=b
\]
of the interval $(a,b)$, such that the complex-valued $V\in
L^{\infty}$ satisfies:
\par
   (i) $\{\Re V(x):x\in(a,b)\}\geq k$ for some $k\in\R$.
\par
   (ii) $V\in C^2$ on each sub-interval $(x_j,x_{j+1})$, $j=0,\ldots,n-1$.
\par
   (iii)
\[
\int_{a}^{b}\abs{\frac{q''(x)}{q^{3/2}(x)}-
\frac{5}{4}\frac{q'^2(x)}{q^{5/2}(x)}}\;\rmd t<\infty
\]
holds, where $q(x):=V(x)-\lam$. Then, for any
$\lam\in\C\backslash\Phi(V)$, we have
\[
\limsup_{h\to 0}\norm(K_h-\lam)^{-1}\norm
\leq\frac{1}{\dist(\lam,\Phi(V))}.
\]
\end{theorem}
\Proof We first prove the case when $-\infty<a<b<+\infty$, noting
that $(iii)$ is then automatically satisfied. Our proof will
involve adding extra points to those in the given partition of the
real interval $(a,b)$. We will then estimate the resolvent norm on
each of the sub-intervals $(x_j,x_{j+1})$. Without loss of
generality therefore, we may initially assume that $V$ is twice
continuously differentiable and bounded on the interval $(-1,1)$,
since the extension to the general case follows easily. Our method
uses the so called WKB approximations \cite{Fed93,Olver74} to the
solutions of the differential equation
\begin{equation}\label{prob1}
(K_h-\lam)f(x)=0
\end{equation}
in the semi-classical limit $h\to 0$. The operator $K_h$ is
defined formally by (\ref{K}). By considering only
$\lam\notin\Phi(V)$ we have $V(x)-\lam\neq 0$ for all
$x\in(-1,1)$. Moreover, since the path $\gam$ defined by
\[
\gam:(-1,1)\to V(x)-\lam
\]
does not cross the negative real axis, condition $(ii)$ also
ensures that we may choose a twice continuously differentiable
branch of $\sqrt{V(x)-\lam}$ such that
\begin{equation}\label{sqrt1}
\Re\sqrt{V(x)-\lam}>0
\end{equation}
for all $x\in(-1,1)$. Denoting the definite integral
\[
\xi(x):=\int_a^x\sqrt{V(t)-\lam}\;\rmd t
\]
where $a\in(-1,1)$ is arbitrary and introduces a constant term
which we will omit from our later calculations, it follows from
(\ref{sqrt1}) that the function $x\mapsto\Re\xi(x)$ is increasing
on the interval $(-1,1)$. This fact will be called upon several
times in our proof. Property (\ref{sqrt1}) will greatly simplify
our application of the WKB approximations, since questions about
the Stokes' phenomenon and valid domains do not then arise.

For fixed $h>0$, let the functions $g_1$ and $g_2$ be linearly
independent, $exact$ classical solutions to (\ref{prob1}). For any
$\alp\in[-1,1]$, put
\[
g\{\alp;x\}:=g_2(\alp)g_1(x)-g_2(x)g_1(\alp)\qquad x\in [-1,1],
\]
so that $g\{-1;x\}$ and $g\{1;x\}$ are also independent
(classical) solutions satisfying $g\{-1;-1\}=g\{1;1\}=0$. Then by
elementary Sturm-Liouville theory, for any $\lam\in\C$ which is
not an eigenvalue, the Green function is given by
\[
G_{\lam}(x,y)=-W_{\lam}^{-1}\left\{\begin{array}{ll}
g\{-1;x\}g\{1;y\} &\mbox{ for $-1\leq x<y$ }\\ g\{1;x\}g\{-1;y\}
&\mbox{ for $y<x\leq 1$ }\\
\end{array} \right.
\]
where the Wronskian
\[
W_{\lam}:=g\{-1;1\}(g_2(0)g_1'(0)-g_1(0)g_2'(0)).
\]
Now consider the operator $\tilde K_h$, again defined formally by
(\ref{K}), but with the extra `boundary condition' $f(0)=0$, so
that $\tilde K_h$ effectively acts on the space $L^2(-1,0)\oplus
L^2(0,1)$. Then the difference resolvent operator
\begin{equation}\label{diffres}
(\tilde K_h-\lam)^{-1}-(K_h-\lam)^{-1}
\end{equation}
is of rank one, provided $\lam$ is not an eigenvalue. Moreover, by
Lemma \ref{2t10} in the Appendix, the operator has a resolvent
kernel given by
\[
\Psi(x,y):=c\phi(x)\phi(y)
\]
where
\[
\phi(x):=\left\{\begin{array}{ll}g\{-1;x\}/g\{-1;0\} &\mbox{ for
$-1\leq x\leq 0$ }\\ g\{1;x\}/g\{1;0\} &\mbox{ for $0\leq x\leq 1$
}\\
\end{array} \right.
\]
and
\[
c:=W_{\lam}^{-1}g\{-1;0\}g\{1;0\}.
\]
Here we have chosen the normalising constant $c$ so that
$\phi(0)=1$. (\ref{diffres}) is seen to be the rank one operator
which acts as
\begin{eqnarray*}
((\tilde K_h-\lam)^{-1}-(K_h-\lam)^{-1})f(x)&:=&
c\phi(x)\int_{-1}^1f(y)\phi(y)\rmd y\\
         &=&c\phi\la f,\overline{\phi}\ra\\
\end{eqnarray*}
on $f\in L^2(-1,1)$, where
\[
\norm(\tilde K_h-\lam)^{-1}-(K_h-\lam)^{-1}\norm
=\abs{c}\norm\phi\norm_2\norm\overline{\phi}\norm_2
=\abs{c}\norm\phi\norm_2^2.
\]
We now turn to the semi-classical behaviour as $h\to 0$. Under the
assumptions of the theorem, asymptotic approximations to solutions
of (\ref{prob1}) are given by (e.g. \cite[p33]{Fed93})
\begin{equation}\label{y1}
y_1(x)=(V(x)-\lam)^{-1/4}\exp\{h^{-1}\xi(x)\}(1+O(h))
\end{equation}
and
\begin{equation}\label{y2}
y_2(x)=(V(x)-\lam)^{-1/4}\exp\{-h^{-1}\xi(x)\}(1+O(h))
\end{equation}
as $h\to 0$. The bound for the remainder term is uniform on
$x\in(-1,1)$, in the sense that $\abs{O(h)}\leq mh$ for $h\leq 1$,
where $m$ does not depend upon $x$. These approximations may be
differentiated with respect to $x$, giving
\begin{equation}\label{y'1}
y_1'(x)=h^{-1}(V(x)-\lam)^{1/4}\exp\{h^{-1}\xi(x)\}(1+O(h))
\end{equation}
and
\begin{equation}\label{y'2}
y_2'(x)=-h^{-1}(V(x)-\lam)^{1/4}\exp\{-h^{-1}\xi(x)\}(1+O(h))
\end{equation}
as $h\to 0$, again the remainder term being uniform on
$x\in(-1,1)$. The functions (\ref{y1}), (\ref{y2}), (\ref{y'1})
and (\ref{y'2}) will be used to estimate the norm of the rank one
difference operator (\ref{diffres}) as $h\to 0$. Indeed,
substituting the approximate solutions $y_1$, $y_2$, $y_1'$ and
$y_2'$ for the exact solutions $g_1$, $g_2$, $g_1'$ and $g_2'$,
one obtains for $\alp,\bet\in(-1,1)$, $\alp<\bet$,
\begin{eqnarray*}
g\{\alp;\bet\}&:=&g_2(\alp)g_1(\bet)-g_2(\bet)g_1(\alp)\\
        &=&y_2(\alp)y_1(\bet)-y_2(\bet)y_1(\alp)\\
        &=&2(V(\alp)-\lam)^{-1/4}(V(\bet)-\lam)^{-1/4}\sinh(h^{-1}
        (\xi(\bet)-\xi(\alp)))
        (1+O(h))\\
        &=&(V(\alp)-\lam)^{-1/4}(V(\bet)-\lam)^{-1/4}\exp\{h^{-1}
        (\xi(\bet)-\xi(\alp))\}
        (1+O(h))\\
\end{eqnarray*}
as $h\to 0$. The last line uses the fact that the function
$x\mapsto\Re\xi(x)$ is increasing on the interval $(-1,1)$, and
that the vanishing term in the sinh function decreases (much) more
rapidly than $O(h)$. The Wronskian simplifies to
\begin{eqnarray*}
W_{\lam}&:=&g\{-1;1\}(g_2(0)g_1'(0)-g_1(0)g_2'(0))\\
        &=&g\{-1;1\}(y_2(0)y_1'(0)-y_1(0)y_2'(0))(1+O(h))\\
        &=&g\{-1;1\}2h^{-1}(1+O(h))\\
\end{eqnarray*}
as $h\to 0$, enabling us to estimate
\begin{eqnarray*}
c&:=&W_{\lam}^{-1}g\{-1;0\}g\{1;0\}\\
 &=&hg\{-1;0\}g\{1;0\}/2g\{-1;1\}(1+O(h))\\
 &=&\frac{h(V(0)-\lam)^{-1/2}\exp\{h^{-1}(\xi(1)-\xi(-1))\}(1+O(h))}
 {2\exp\{h^{-1}(\xi(1)-\xi(-1))\}(1+O(h))}\\
 &=&O(h)\\
\end{eqnarray*}
as $h\to 0$. It also follows that for $0< x\leq 1$
\begin{eqnarray*}
\phi(x)&:=&g\{1;x\}/g\{1;0\}\\
       &=&\frac{(V(1)-\lam)^{-1/4}(V(x)-\lam)^{-1/4}\exp\{h^{-1}(\xi(1)-\xi(x))\}(1+O(h))}
       {(V(1)-\lam)^{-1/4}(V(0)-\lam)^{-1/4}\exp\{h^{-1}(\xi(1)-\xi(0))\}(1+O(h))}\\
       &=&\frac{(V(x)-\lam)^{-1/4}}{(V(0)-\lam)^{-1/4}}\exp\{h^{-1}
       (\xi(0)-\xi(x))\}(1+O(h))\\
\end{eqnarray*}
as $h\to 0$. We note that
\[
\Re(\xi(0)-\xi(x))<0
\]
for $0<x\leq 1$; and $\phi(0)=1+O(h)$. Then, using the fact that
$\phi$ is even, and applying the method of steepest descents
\begin{eqnarray*}
\norm\phi\norm_{L^2(-1,1)}^2&=& 2\int_0^1\abs{\phi(x)}^2\rmd x\\
 &\leq&
 c_1\int_0^1\abs{\exp\{h^{-1}(\xi(0)-\xi(x))\}(1+O(h))}^2\rmd x\\
 &=&c_1\int_0^1\exp\{2h^{-1}\Re(\xi(0)-\xi(x))\}\rmd x\;(1+O(h))\\
 &=&c_1\int_0^1\exp\{-2h^{-1}\Re(\xi'(0)x))\}\rmd x\;(1+O(h))\\
 &=&\frac{c_1 h}{2\Re(\xi'(0))}(1+O(h))\\
 &=&O(h)\\
\end{eqnarray*}
as $h\to 0$, since $\Re(\xi'(0))=\Re\sqrt{V(0)-\lam}>0$. Therefore
\begin{eqnarray*}
\norm(\tilde K_h-\lam)^{-1}-(K_h-\lam)^{-1}\norm
&=&\abs{c}\norm\phi\norm_2^2\\ &=&O(h)\norm\phi\norm_2^2\\
&=&O(h^2)\\
\end{eqnarray*}
as $h\to 0$.

Now defining the operator $\tilde K_h$ formally by (\ref{K}) but
with the finite number of boundary conditions
\[
f(x_j)=0\qquad j=0,\ldots,n,
\]
$\tilde K_h$ then effectively acts on
\[ L^2(-1,x_1)\oplus L^2(x_1,x_2)\oplus\cdots\oplus
L^2(x_{n-1},1).
\]
By induction, our argument so far shows that the operator
\[
(\tilde K_h-\lam)^{-1}-(K_h-\lam)^{-1}
\]
is of at most rank $n$, and we have the norm resolvent convergence
\begin{equation}\label{norms}
\lim_{h\to 0}\norm(\tilde K_h-\lam)^{-1}-(K_h-\lam)^{-1}\norm=0.
\end{equation}
Moreover, letting $V_j$ denote the potential $V$ restricted to the
interval $(x_j,x_{j+1})$, condition (i) ensures that one can apply
Proposition \ref{p1} to $\tilde K_h$ separately on each interval
$(x_j,x_{j+1})$. Therefore, taking
\[
\lam\in\C\backslash\bigcup_{j=0}^{n-1}\conv(\Phi(V_j))
\]
we have
\begin{equation}\label{bounds}
\norm(\tilde K_{h}-\lam)^{-1}\norm \leq\max_{j}
\left\{\frac{1}{\dist(z,\conv(\Phi(V)))}\right\}
<\infty
\end{equation}
uniformly on $h>0$. As the partition of $(-1,1)$ becomes
increasingly fine, it is clear that for every
$\lam\in\C\backslash\Phi(V)$ one has
\[
\min_{j}\{\dist(\lam,\conv(\Phi(V_j)))\}\longrightarrow
\dist(\lam,\Phi(V)).
\]
Then, by (\ref{norms}) and (\ref{bounds}) it follows that
\[
\limsup_{h\to 0}\norm(K_{h}-\lam)^{-1}\norm \leq
\frac{1}{\dist(\lam,\Phi(V))}
\]
to complete the proof in the finite interval case.

On the interval $(a,+\infty)$, the WKB approximations to the
classical solutions of (\ref{prob1}) take the form
\[
\tilde y_1(x)=(V(x)-\lam)^{-1/4}\exp\{h^{-1}\xi(x)\}(1+o(1))
\]
and
\[
\tilde y_2(x)=(V(x)-\lam)^{-1/4}\exp\{-h^{-1}\xi(x)\}(1+o(1))
\]
as $x\to\infty$, for all $h>0$, provided $(iii)$ holds (see
\cite[p50]{Fed93}). Comparing $\tilde y_{1,2}$ with $y_{1,2}$ in
the proof above, and noting that $\tilde y_2(x)\to 0$
exponentially as $x\to +\infty$, one can check that the proof just
given still carries through. The interval $(-\infty,b)$ is dealt
with similarly, by a change of signs, completing the proof for the
general case.


\section{A Different Approach}

In Theorem \ref{t6} we relied upon the theory of ODEs to prove the
norm resolvent convergence (\ref{norms}). The proof of the next
theorem uses a powerful but technically simple construction, the
so-called `Twisting Trick' \cite[Section 8.6]{Dav95} or
\cite{Dav82}.

\begin{theorem}\label{t7}
Let $K_h$ be defined by (\ref{K}) on $L^2(\Omega)$, where $\Omega$
is a bounded region in $\R^N$ and $V:\bar\Omega\to\C$ is
continuous. Then, for $\lam\notin\Phi(V)$
\[
\limsup_{h\to 0}\norm(K_h-\lam)^{-1}\norm
\leq\frac{1}{\dist(\lam,\Phi(V))}.
\]
\end{theorem}

\Proof If $\lam\notin\conv(\Phi(V))$ then Proposition \ref{p1}
applies and we are done. So, assume that
$\lam\in\conv(\Phi(V))\backslash\Phi(V)$ is given. For any
$\del>0$ we define $\{S_j\}$ to comprise $N$-dimensional cubes of
the form $\{(x_1,\ldots,x_N):\delta r_i< x_i<\delta(r_i+1)\}$,
where the $r_1,\ldots,r_N$ take integer values. Then by the
uniform continuity of $V$ on $\bar\Omega$, there exists a covering
of $\bar\Omega$ by disjoint cubes
\[
\bar\Omega\subseteq\bigcup_{j=1}^M\bar S_j
\]
each of side length $\del$, such that
\[
\lam\notin\bigcup_{j=1}^M\conv\left\{\Phi\left(V\big|_{\bar\Omega\cap\bar
S_j} \right)\right\}.
\]
In addition, for any given $0<\alp<1$ we can always take
$\delta>0$ small enough so that
\begin{equation}\label{alpest}
\dist\left(\lam,\bigcup_{j=1}^M\conv\left\{\Phi\left(V\big|_{\bar\Omega\cap\bar
S_j} \right)\right\}\right)\geq\alpha\dist(\lam,\Phi(V)).
\end{equation}
The proof proceeds by a series of bisections in each of the $N$
dimensions of $\Omega$. Choose a point $c=(c_1,\ldots,c_N)$ such
that each $c_i=\del r_i$ for some integer $r_i$. Then the
hyperplane $\{x:x_i=c_i\}$ splits the cubes into two families; one
covering the region $\Omega_1:=\{x:x_i-c_i>0\}$, and the other
covering the region $\Omega_2:=\{x:x_i-c_i<0\}$. Define
\[
V_1(x):=\left\{\begin{array}{ll} V(x) &\mbox{ if $x\in\Omega_1$ }
\\ m &\mbox{ if $x\in\Omega_2$ }\\
\end{array} \right.
\]
and
\[
V_2(x):=\left\{\begin{array}{ll} V(x) &\mbox{ if $x\in\Omega_2$ }
\\ m &\mbox{ if $x\in\Omega_1$ }\\
\end{array} \right.
\]
where $m$ is some sufficiently large real number. Then consider
the two operators
\[
H_1:=\pmatrix{ -h^2\lap +V & 0 \cr 0 & -h^2\lap+m \cr },
H_2:=\pmatrix{ -h^2\lap +V_1 & 0 \cr 0 & -h^2\lap +V_2 \cr },
\]
both acting in the Hilbert space $\cH:=L^2(\Omega)\oplus
L^2(\Omega)$. We will show that
\[
\norm(H_1-\lam)^{-1}\norm -\norm(H_2-\lam)^{-1}\norm\to 0
\]
as $h\to 0$. Defining $\theta:\R\to[0,\pi/2]$ by
\[
\theta(s):=\left\{\begin{array}{ll} \pi/2 &\mbox{ if $s\leq -1/3$
} \\ \pi(1-3s)/4 &\mbox{ if $-1/3\leq s\leq 1/3$ }
\\ 0 &\mbox{ if $s\geq 1/3$ }
\\
\end{array} \right.
\]
we define the unitary operator $U_h:\cH\to\cH$ by
\[
U_h\pmatrix { f(x) \cr g(x) \cr }:= \pmatrix {
\cos\theta((x_i-c_i)/h^{\gam}) & \sin\theta((x_i-c_i)/h^{\gam})
\cr -\sin\theta((x_i-c_i)/h^{\gam}) &
\cos\theta((x_i-c_i)/h^{\gam}) \cr }\pmatrix { f(x) \cr g(x) \cr }
\]
where $\gam>0$ is to be determined. Thus
\[
U_h(x)=\left\{\begin{array}{ll} \pmatrix { 1 & 0 \cr 0 & 1 \cr }
&\mbox{ if $x_i\geq c_i+h^{\gam}/3$ }\\ \pmatrix { 0 & 1 \cr -1 &
0 \cr } &\mbox{ if $x_i\leq c_i-h^{\gam}/3$ }.\\
\end{array} \right.
\]
To ease notation we denote the following functions, which are to
be regarded as multiplication operators on $L^2(\Omega)$,
\[
C_h(x):=\cos\theta((x_i-c_i)/h^{\gam})
\]
\[
S_h(x):=\sin\theta((x_i-c_i)/h^{\gam})
\]
together with the partial differentiation operator $D_i:=\partial
/\partial x_i$. Then, one can show using elementary matrix
calculations that
\begin{equation}\label{8.6.1}
U_hH_1U_h^*=H_2+P_hD_i+Q_h+G_h
\end{equation}
where $P_h$, $Q_h$ and $G_h$ are the matrix-valued functions on
$\Omega$ given by
\[
P_h:=\frac{3\pi h^{2-\gam}}{2}\chi_{\cC}\pmatrix{ 0 & 1 \cr -1 & 0
\cr }
\]

\[
Q_h:=\frac{9\pi^2h^{2-2\gam}}{16}\chi_{\cC}\pmatrix{ 1 & 0 \cr 0 &
1 \cr }
\]
and
\[
G_h:=(m-V)\pmatrix{ \chi_{\Omega_1}S_h^2-\chi_{\Omega_2}C_h^2 &
C_hS_h \cr C_hS_h & -\chi_{\Omega_1}S_h^2+\chi_{\Omega_2}C_h^2 \cr
},
\]
where
\[
\cC:=\{x\in\Omega:-h^{\gam}/3<x_i-c_i<h^{\gam}/3\}.
\]
Thus $\norm P_h\norm =O(h^{2-\gam})$, $\norm Q_h\norm
=O(h^{2-2\gam})$ and since $\chi_{\Omega_1}S_h$ and
$\chi_{\Omega_2}C_h$ are $O(h^{\gam})$, $\norm G_h\norm
=O(h^{\gam})$, as $h\to 0$. We therefore take the optimal value
$\gam=2/3$. Now, for our given $\lam$, one may write
\[
(H_2-\lam)^{-1}-(U_hH_1U_h^*-\lam)^{-1}=(H_2-\lam)^{-1}(P_hD_i+Q_h+G_h)
(U_hH_1U_h^*-\lam)^{-1},
\]
so that
\begin{equation}\label{8.6.3}
\norm(H_2-\lam)^{-1}-(U_hH_1U_h^*-\lam)^{-1}\norm\leq\norm(H_2-\lam)^{-1}\norm
\norm G_h\norm+\norm P_h\norm \norm
D_i(U_hH_1U_h^*-\lam)^{-1}\norm +\norm Q_h\norm.
\end{equation}
Now $(H_1-\lam)^{-1}$ and $(H_2-\lam)^{-1}$ are bounded from $\cH$
to $W_0^{1,2}$, $D_i$ is bounded from $W_0^{1,2}$ to $\cH$, and
$U_h$ is uniformly bounded from $W_0^{1,2}$ to $W_0^{1,2}$ for all
$h\leq 1$. Thus
\begin{eqnarray*}
\norm D_i(U_hH_1U_h^*-\lam)^{-1}\norm&=&\norm
D_iU_h(H_1-\lam)^{-1}U_h^*\norm\\ &\leq &\bet
\end{eqnarray*}
for some $\bet<\infty$ and all $h\leq 1$. Therefore, from
(\ref{8.6.3}), we obtain
\[
\norm(U_hH_1U_h^*-\lam)^{-1}-(H_2-\lam)^{-1}\norm=O(h^{1/2})
\]
as $h\to 0$, so that
\[
\norm(H_1-\lam)^{-1}\norm=\norm(U_hH_1U_h^*-\lam)^{-1}\norm
                         =\norm(H_2-\lam)^{-1}\norm +O(h^{1/2})
\]
where, applying Proposition \ref{p1}, one has
\[
\norm(H_2-\lam)^{-1}\norm
\leq\max\left\{\frac{1}{\dist(\lam,\conv(\Phi(V_1)))},\frac{1}{
\dist(\lam,\conv(\Phi(V_2)))}\right\}.
\]
One can now bisect each of $\Omega_1$ and $\Omega_2$ in the same
manner, and carry out the above process on four copies of
$L^2(\Omega)$. Repeating until $\Omega$ has been divided into
`strips' of thickness $\del$, one changes to another coordinate
direction and repeats the process until all $N$ dimensions have
been decomposed. Then, recalling (\ref{alpest}) we have
\[
\norm(H_h-\lam)^{-1}\norm\leq\frac{\alp^{-1}}{\dist(\lam,\Phi(V))}+O(h^{1/2})
\]
as $h\to 0$, where $0<\alp<1$ was arbitrarily chosen. Taking
$\alp$ as close as one likes to $1$ will then complete the proof.


\section{Appendix}

We give a proof of the following well-known result.

\begin{lemma}\label{2t10}
Let $L$ be the Sturm-Liouville operator
\begin{equation}\label{Ldefn}
(L-\lam)f(x):=-\frac{\rmd^2f(x)}{\rmd x^2}+V(x)f(x)-\lam f(x)=0
\end{equation}
acting in $L^2(a,b)$, together with boundary conditions
\begin{equation}\label{bdrycon}
f(a)=f(b)=0.
\end{equation}
Here V is a complex-valued continuous function on $[a,b]$, and
$\lam$ is a complex constant. We may assume that $a<0<b$, and let
$\tilde L$ denote the operator given formally by (\ref{Ldefn}) but
subject to the additional condition $f(0)=0$. Then
\[
(\tilde L-\lam)^{-1}-(L-\lam)^{-1}
\]
is a rank one operator for all $\lam\notin\{\Spec(\tilde
L)\cup\Spec(L)\}$.
\end{lemma}
\Proof Let $u$, $v$ be a pair of independent classical solutions
to the differential equation (\ref{Ldefn}), with $u$, $v$
satisfying the boundary conditions at $a$, $b$ respectively.
Multiplying $u$ by the constant $v(0)/u(0)$, we may further assume
that $u(0)=v(0)$. Then the Green function is given by
\[
G_{\lam}(x,y):=-W_{\lam}^{-1}\left\{\begin{array}{ll} u(x)v(y)
&\mbox{ if $a\leq x\leq y\leq b$} \\ v(x)u(y) &\mbox{ if $b\geq
x\geq y\geq a.$}
\end{array}\right.
\]
For any $x\in(a,b)$, the Wronskian is given by
\[
W_{\lam}:=u(x)v'(x)-u'(x)v(x).
\]
Imposing the boundary condition $f(0)=0$, $\tilde L$ effectively
acts on the Hilbert space $L^2(a,0)\oplus L^2(0,b)$. Then, putting
\[
w(x):=u(x)-v(x)
\]
we see that $w(0)=0$, and the new Wronskian on $L^2(a,0)$:
\begin{eqnarray*}
W^-&:=&u(0)w'(0)-u'(0)w(0) \\
   &=&u(0)(u'(0)-v'(0))-u'(0)(u(0)-v(0)) \\
   &=&-u(0)v'(0)+u'(0)v(0) \\
   &:=&-W.
\end{eqnarray*}
By a similar calculation, the Wronskian on $L^2(0,b)$ is given by
$W^+=W$. Hence we construct the new `Green' function
\[
\tilde k(x,y):=-W^{-1}\left\{\begin{array}{ll} 0 &\mbox{ if $a\leq
x\leq 0$ and $b\geq y\geq 0$}\\ 0 &\mbox{ if $b\geq x\geq 0$ and
$a\leq y\leq 0$}
\\ w(y)v(x) &\mbox{ if $0\leq y\leq x\leq b$}\\
v(y)w(x) &\mbox{ if $0\leq x\leq y$}\\ -w(x)u(y) &\mbox{ if $a\leq
y\leq x\leq 0$}\\ -u(x)w(y) &\mbox{ if $a\leq x\leq y\leq 0.$}
\end{array}\right.
\]
It is then straightforward to check that putting
\[
\sig(x):=\left\{\begin{array}{ll} v(x) &\mbox{ if $b\geq x\geq 0
$}\\ u(x)&\mbox{ if $a\leq x\leq 0$}
\end{array}\right.
\]
the `Green' function of
\[
(\tilde L-\lam)^{-1}-(L-\lam)^{-1}
\]
is given by
\[
W^{-1}\sig(x)\sig(y)
\]
for all $x,y$ in $(a,b)$. Clearly $\sig\in L^2(a,b)$, and for all
$f\in L^2(a,b)$
\begin{eqnarray*}
\left((\tilde L-\lam)^{-1}-(L-\lam)^{-1}\right)f(x)&=&
W^{-1}\int_a^b\sig(x)\sig(y)f(y)\;\rmd y\\
&=&W^{-1}\sig(x)\int_a^b\sig(y)f(y)\;\rmd y.
\end{eqnarray*}
\[
(\tilde L-\lam)^{-1}-(L-\lam)^{-1}
\]
is a compact rank one operator on $L^2(a,b)$, which completes the
proof.

\vskip 0.5in
{\bf Acknowledgements }I should like to thank A. Aslanyan, E. B.
Davies and Y. Safarov for valuable discussions and comments.
\par

\vskip 0.3in Department of Mathematics \newline King's College
\newline Strand \newline London WC2R 2LS \newline England \\
e-mail:Redparth@mth.kcl.ac.uk \vfil
\end{document}